\documentclass[12pt]{article}
\usepackage[utf8]{inputenc}
\usepackage[english]{babel}
\usepackage{amsfonts}
\usepackage{amssymb}
\usepackage{amsmath}
\usepackage{amsthm}
 \usepackage{CJK}
\numberwithin{equation}{section}
\begin{document}

 \title{Jacobian of incidence schemes }
\author{ B. Wang \\
\begin{CJK}{UTF8}{gbsn}
(汪      镔)
\end{CJK}}
\date{Aug 9, 2017}

\newcommand{\Addresses}{{
  \bigskip
  \footnotesize

   \textsc{Department of Mathematics, Rhode Island college, Providence, 
   RI 02908}\par\nopagebreak
  \text{E-mail address}:  \texttt{binwang64319@gmail.com}

}}

\maketitle

\newtheorem{thm}{Theorem}[section]

\newtheorem{proposition}[thm]{\bf {Proposition} }
\newtheorem{theorem}[thm]{\bf {Theorem} }
\newtheorem{ex}[thm]{\bf Example}
\newtheorem{corollary}[thm]{\bf Corollary}
\newtheorem{definition}[thm]{\bf Definition}
\newtheorem{lemma}[thm]{\bf Lemma}

\newcommand{\xdownarrow}[1]{%
  {\left\downarrow\vbox to #1{}\right.\kern-\nulldelimiterspace}
}

\newcommand{\xuparrow}[1]{%
  {\left\uparrow\vbox to #1{}\right.\kern-\nulldelimiterspace}
}

\begin{abstract}  
We show a direct calculation of Jacobian matrices  in the old problems of 
rational curves on generic hypersurfaces.
\end{abstract}

\tableofcontents

\section{Introduction}

Let's consider the problems in $\mathbf P^4$, the $\mathbb C$-projective space of dimension $4$.

\bigskip

Old problems: \par
(1) Is there a rational curve on a generic quintic hypersurface $X \subset \mathbf P^4$?\par
(2) If the answer is ``yes", how do the rational curves lie on $X$? 

\bigskip

Two questions have been answered by H. Clemens and S. Katz ([1], [2]).

\begin{theorem}  (Clemens and Katz)
There are smooth rational curves $C$ of each degree $d>0$ on a generic quintic hypersurface $X$ and
\begin{equation}
N_{C/X}\simeq \mathcal O_{\mathbf P^1}(-1)\oplus  \mathcal O_{\mathbf P^1}(-1).
\end{equation}
\end{theorem}

Let's reformulate the theorem.
Let \begin{equation}M=( H^0(\mathcal O_{\mathbf P^1}(d))^{\oplus 5}\simeq \mathbb C^{5d+5},\end{equation}
the parameter space of parametrized rational curves, i.e. the space of all
5 tuples of polynomials in one variable of degree $d$.  
 An open set $M_d$ of $M$ corresponds (not equal) to the set \begin{equation} 
\{ \ maps\ c: \mathbf P^1\to \mathbf P^4: embedding\ to \ its \ image, \ deg(c(\mathbb P^1))=d \}.\end{equation}
Let $X$ be  a quintic 3-fold in $\mathbb P^4$, i.e 
\begin{equation}
X=div(f), f\in H^5(\mathcal O_{\mathbf P^4}(5)).
\end{equation}

Let the  incidence scheme be 
$$ I_{X}=\{ \ maps\ c\in M_d: c(\mathbf P^1)\subset X \}\subset M_d. $$ 
Notice $I_X$ is well-defined and its defining equations are parametrized algebraically by 
$f\in  H^5(\mathcal O_{\mathbf P^4}(5) $ with $div(f)=X$. Hence its
Jacobian  matrix $J_X$  of defining equations  is also well-defined, and
each entry is a polynomial in variables $f\in H^5(\mathcal O_{\mathbf P^4}(5)$ and $c\in M$.
The Jacobian matrix $J_X$ is not an invariant of $I_X$. In the following we give 
 a definition used in this paper.  Notice 
$f(c(t))$ is a polynomial of degree $5d$ in $t$, i.e.  
\begin{equation}
\sum_{j=0}^{5d}k_j(c, f) t^j=f(c(t).
\end{equation}
Hence $ I_{X}$ is the projection of
$$\{ (c, f)\in M_d\times  H^5(\mathcal O_{\mathbf P^4}(5)): k_j(c, f)=0, j=0, \cdots, 5d\}$$ to
$M_d$. 

\bigskip

\begin{definition} \quad\par
(1) Let $\theta_1, \cdots, \theta_{5d+5}$ be affine coordinates for $M$.
Define 
\begin{equation}
J_{X}={\partial (k_0(c, f), \cdots, k_{5d}(c, f))\over \partial (\theta_1, \cdots, \theta_{5d+5})}
\end{equation}

(2) More generally, let $k_1, \cdots, k_m$ be  polynomials in $M$ defining $I_X$. 
We use the same notation $J_X$ to denote the Jacobian matrix

\begin{equation}
{\partial (k_1, \cdots, k_m)\over \partial (\theta_1, \cdots, \theta_{5d+5})}
\end{equation}

\end{definition}

\bigskip

Thus the matrix $J_{X}$ is not well-defined as an invariant of $I_X$. 
But its rank at $c\in I_X$ is  well-defined intrinsically as 
\begin{equation} 5d+5-dim(T_c I_X).\end{equation}
We say the matrix is non-degenerate if it has the full rank.
\par

Then theorem 1.1 says
\bigskip

\begin{proposition}\quad   Let  $X$ be  generic in $ \mathbf P(H^5(\mathcal O_{\mathbf P^4}(5))$. \par
There is a $c\in I_X$ such that  $J_X$ is non-degenerate at $c$.
\end{proposition}

\bigskip
To construct such a $c$, Clemens considered a special type of quintics $X_0$, each of which contains  a smooth K3 surface $Q$ of
$\mathbf P^3$. Then he proceeded to prove 
\bigskip

(1) For each of infinitely many $d$, $Q$ contains a smooth rational curve \par\hspace {1cc}  $c_0\in M_d$;\par
(2) Applying the structure of normal bundles of $c_0(\mathbf P^1)$, he obtained that   \par\hspace {1cc} $J_{X_0}$ is non-degenerate at $c_0$;\par
(3) The non-degeneracy  is extended to generic $X$. 
\bigskip

S. Katz later added a new result of Mori :  for each degree $d$, a smooth K3 surface containing a smooth rational curve of degree $d$ always exists. This extends Clemens' result to all degrees $d$.

\bigskip

In the  following we give a new  approach to (2):  A DIRECT CALCULATION of $J_{X_0}$ at $c_0$ without using the normal bundles. 

\bigskip

\section{New approach to Jacobian matrices}

\bigskip

We present the direct calculation. \bigskip

\begin{proposition}
Let $Q$ be a smooth quartic hypersurface in $\mathbf P^3\subset \mathbf P^4$. Assume there is  $c_0\in M_d$ on $Q$. Then there is a special quintic hypersurface $X_0\subset \mathbf P^4$ containing $c_0(\mathbf P^1)$, such that
the Jacobian $J_{X_0}$ is non-degenerate at $c_0$.

\end{proposition}

\bigskip

\begin{proof}

We use affine coordinate $t\in \mathbb C\subset \mathbf P^1$ for $\mathbf P^1$.\par

Construction of $f_0$ ( $div(f_0)=X_0$): 
Let $z_0, \cdots, z_4$ be homogeneous coordinates of $\mathbf P^4$, $q(z)$ a smooth homogeneous quartic in $z_0, \cdots, z_3$ 
and $z_4=0=q(z)$ is $Q$. Let 
$p(z)$ a generic homogeneous quartic in $z_0, \cdots, z_3, z_4$, and  $l(z)$ be a generic, homogeneous  linear polynomial in $z_0, \cdots, z_4$.  
Let 
\begin{equation}
f_0=l(z) q(z)+z_4 p(z)
\end{equation}
be the quintic containing $c_0$. 

\bigskip

Calculation of Jacobian $J_{f_0}$:   Let \begin{equation}
y_4^0, \cdots, y_4^d, y_1, \cdots, y_{4d+4} 
\end{equation}
be coordinates of $M$ around the point $c_0$ where $\sum_{i=0}^dy_4^i t^i$ is the polynomial for
 the last copy $H^0(\mathcal O_{\mathbf P^1}(d))$, and 
$y_i$ are for the rest 4 copies of $H^0(\mathcal O_{\mathbf P^1}(d))$.
Next we explain the defining polynomials of $I_{X_0}$ on $M$.
Choose $5d+1$ distinct points 
$$t_1, \cdots, t_{5d+1}$$ 
in $\mathbf P^1$ in the following way:
 Let $t_1, \cdots, t_{d}$ be distinct $d$ zeros of $l(c_0(t))=0$, but not zeros of $p(c_0(t))=0$,  $(t_{d+1}, t_{d+2}, \cdots, t_{5d+1})$  
are generic in $\Pi_{4d}\mathbf P^1$. Using these points, 
 we obtain $5d+1$ homogeneous polynomials in $M$, 
\begin{equation}
f_0(c(t_1)),\cdots, f_0(c(t_{5d+1}))
\end{equation}
which are local defining polynomials of $I_{X_0}$ around the point $c_0$.
Then a direct calculation of the Jacobian divides
$J_{X_0}|_{c_0}$  into blocks

\begin{equation}
J_{X_0}=\left(  \begin{array}{cc} \mathcal A_{11} & 0  \\
\mathcal A_{21}   &\mathcal A_{22}\\
\end{array}\right)
\end{equation}
where the 
\begin{equation}
\mathcal A_{11}={\partial ( f(c(t_1), \cdots, f(c(t_{d+1}))\over \partial (y_4^0, \cdots, y_4^d)}|_{c_0};\end{equation}

\begin{equation}
\mathcal A_{12}={\partial ( f(c(t_1), \cdots, f(c(t_{d+1}))\over \partial (y_1 \cdots, y_{4d+4})}|_{c_0}=0;\end{equation}

\begin{equation}
\mathcal A_{21}={\partial ( f(c(t_{d+2}), \cdots, f(c(t_{5d+1}))\over \partial (y_4^0, \cdots, y_4^d)}|_{c_0};\end{equation}

\begin{equation}
\mathcal A_{22}={\partial ( f(c(t_2), \cdots, f(c(t_{5d+1}))\over \partial (y_1, \cdots, y_{4d+4})}|_{c_0};\end{equation}

Then taking the partial derivatives, we obtain

\begin{equation}
\mathcal A_{11}=\left(  \begin{array}{cccc} t_1^d p(c_0(t_1)) & \cdots & t_1p(c_0(t_1)) & p(c_0(t_1))\\
\vdots &\cdots & \vdots & \vdots \\
t_{d+1}^d p(c_0(t_{d+1})) & \cdots & t_{d+1}p(c_0(t_{d+1})) & p(c_0(t_{d+1})) \\
\end{array}\right)\end{equation}

\begin{equation}
\mathcal A_{22}=\left(  \begin{array}{ccc}  l(c_0(t_{d+2})) {\partial q\over \partial y_1}|_{c_0(t_{d+2})}
 &\cdots &  l(c_0(t_{d+2})) {\partial q\over \partial y_{4d+4}}|_{c_0(t_{d+2})} \\
\vdots & \cdots  &\vdots\\
 l(c_0(t_{5d+1})) {\partial q\over \partial y_1}|_{c_0(t_{5d+1})}
 &\cdots &  l(c_0(t_{5d+1})) {\partial q\over \partial y_{4d+4}}|_{c_0(t_{5d+1})} \\
\end{array}\right)\end{equation}

$\mathcal A_{11}$ is obviously invertible (because $p(c_0(t_i))\neq 0 , 1\leq i\leq d$). Thus it suffices to
show $\mathcal A_{22}$ has full rank. Furthermore it is sufficient to show
\begin{equation}
\mathcal A_0=\left(  \begin{array}{ccc} {\partial q\over \partial y_1}|_{c_0(t_{d+2})} &\cdots
 & {\partial q\over \partial y_{4d+4}}|_{c_0(t_{d+2})} \\
\vdots & \cdots  &\vdots\\
{\partial q\over \partial y_1}|_{c_0(t_{5d})} &\cdots
 & {\partial q\over \partial y_{4d+4}}|_{c_0(t_{5d})}\\
 {\partial q\over \partial y_1}|_{c_0(t_{5d+1})} &\cdots
 & {\partial q\over \partial y_{4d+4}}|_{c_0(t_{5d+1})}\end{array}\right),\end{equation}
has full rank.
Suppose $\mathcal A_0$ does not have full rank.  Then the solution set 
 to the system of linear equations 
\begin{equation}
\left(  \begin{array}{ccc} {\partial q\over \partial y_1}|_{c_0(t_{d+2})} &\cdots
 & {\partial q\over \partial y_{4d+4}}|_{c_0(t_{d+2})} \\
\vdots & \cdots  &\vdots\\
{\partial q\over \partial y_1}|_{c_0(t_{5d})} &\cdots
 & {\partial q\over \partial y_{4d+4}}|_{c_0(t_{5d})}\\
{\partial q\over \partial y_1}|_{c_0(t_{5d+1})} &\cdots
 & {\partial q\over \partial y_{4d+4}}|_{c_0(t_{5d+1})} 
\end{array}\right)\cdot \left(  \begin{array}{c} x_1 \\
\vdots \\
x_{4d+3}\\
x_{4d+4}
\end{array}\right)
=0\end{equation}
is at least $5$ dimensional. On the other hand, due to the genericity of $t_{5d+1}$, 
each non-zero solution of (2.12) gives a non zero section 
of $T_{Q}|_{C_0}$, where $C_0=c_0(\mathbf P^1)$.  
Therefore  $H^0(T_{Q}|_{C_0})$ is at least $5-1=4$-dimensional ( projetivize $M$). 
But due to the Calabi-Yau nature of $Q$, we know  $H^0(T_{Q}|_{C_0})$ must be isomorphic to
$H^0(T_{C_0})$ which is only 3-dimensional. The contradiction shows $J_{X_0}$ has full rank. 
\end{proof}

\bigskip

\section{ New direction}
This calculation gives a new direction in attacking the Jacobian matrix $J_{X}$:  break the Jacobian matrix into manageable blocks. 
After Clemens', there is a need of new setting to break the Jacobian.  Let's have a new set-up.\bigskip

Let $S=\mathbf P(\mathcal O_{\mathbf P^4}(5)$ be the parameter space of quintic 3-folds in $\mathbf P^4$. 
Thus $X_0$ is viewed as a point  in $S$, i.e. any point $s\in S$ corresponds to a quintic 3-fold $X_s\subset \mathbf P^4$. 
We also extend $M_d$ to all the irreducible rational curves. 
Next we take an open set of plane ${\mathbb L}\subset S$ containing the point $X_0$.
 Then we repeat the definition of a incidence scheme above to define the  secondary incidence scheme $I_{\mathbb L}$ to be an irreducible 
component of 
\begin{equation}
\{  c\in M_d: c( P^1)\subset  X_l, l\in \mathbb L\}\subset M_d. \end{equation}

Then a direct calculation can be applied for the following proposition: 

\bigskip

\begin{proposition} \par Let $\mathbb L$ be generic in $S$. \par
  If for every $l\in \mathbb L$, there is a  $c\in M_d$ such that $$ c(\mathbf P^1)\subset  X_l, $$
\hspace{1cc} then the Jacobian matrix $J_{\mathbb L}$ of $I_{\mathbb L}$ at generic point is non-degenerate.  

\end{proposition}
\bigskip

The direct calculation of the Jacobian  holds for rational curves and generic hypersurfaces in all $\mathbf P^n, n\geq 4$.
Hence the Calabi-Yau condition does not play a role in the calculation.  The proof was contained in another paper.

\end{document}